\documentclass[adm]{birkmult}

\usepackage{amssymb}
\newtheorem{Proposition}{Proposition}[section]

\newtheorem{Lemma}[Proposition]{Lemma}
\newtheorem{Theorem}[Proposition]{Theorem}

\newcommand{\R}{\mathbb{R}}

\DeclareMathOperator{\inte}{int}
\DeclareMathOperator{\grad}{grad}

\renewcommand{\labelenumi}{\roman{enumi})}

\title{Hilbert geometry of polytopes}

\author{Andreas Bernig}
 
\email{andreas.bernig@unifr.ch}

\address{D\'epartement de Math\'ematiques, Chemin du Mus\'ee 23, 1700 Fribourg, Switzerland}

\begin{document}

\begin{abstract}
It is shown that the Hilbert metric on the interior of a convex polytope is bilipschitz to a normed vector space of the same dimension. 
\end{abstract}

\thanks{
\\ Supported
  by the Schweizerischer Nationalfonds grants SNF PP002-114715/1 and 200020-121506/1.}

\maketitle 
\section{Introduction}

Given a compact convex set $K$ in a finite-dimensional vector space $V$, the Hilbert metric (also called Hilbert geometry) on $\inte K$ is defined by 
\begin{displaymath}
 d(x,y):=\frac12 \left|\log [a_1,a_2,x,y]\right|, \quad x,y \in \inte P.
\end{displaymath}
Here $a_1,a_2$ are the intersections of the line through $x$ and $y$ with the boundary of $K$ (if $x=y$, one sets $d(x,y):=0$). Hilbert metrics are examples of projective Finsler metrics, i.e. Finsler metrics such that straight lines are geodesics. Hilbert's fourth problem was to classify projective Finsler metrics. This problem was solved by Pogorelov \cite{pogorelov73}, see also \cite{alvarez05} for a symplectic approach. 

In the last few years, there has been a renewed interest in Hilbert geometries and many research papers were published, see \cite{karlsson_guennadi02,sociemethou04, foertsch_karlsson05,colbois_vernicos07, yvessurvey, vernicos05,guichard05,berck_bernig_vernicos} to cite just a few of them.

A natural question is to classify Hilbert metrics up to bilipschitz maps or up to quasi-isometries. In this direction, it was shown by Colbois and Verovic \cite{colbois_verovic04} that if $\partial K$ is $C^2$ with positive Gauss curvature, then the Hilbert metric is bilipschitz to the $n$-dimensional hyperbolic metric. In a similar spirit, it was shown in \cite{berck_bernig_vernicos} that the volume entropy of an $n$-dimensional convex body with $C^{1,1}$ boundary equals the volume entropy of hyperbolic space, which is $n-1$. Since the volume entropy is an invariant under bilipschitz maps, the Hilbert metric of a body with $C^{1,1}$ boundary can not be bilipschitz to a normed vector space. 

On the other extreme, Hilbert metrics of polygons are bilipschitz to normed spaces, as was shown by Colbois, Vernicos and Verovic \cite{colbois_vernicos_verovic08}. Elaborating an argument of Foertsch and Karlsson \cite{foertsch_karlsson05}, Colbois and Verovic \cite{colbois_verovic08} showed that if the Hilbert metric of a compact convex body is quasi-isometric to a normed vector space (in particular, if it is bilip\-schitz), then the body is a polytope (recall that a polytope is the convex hull of a finite number of points). 

This raised the question whether the converse holds in general, i.e. if the Hilbert metric of any polytope is quasi-isometric (or bilipschitz) to a normed vector space. The aim of this note is to answer this question in the positive. 

\begin{Theorem} \label{thm_main_thm_exact}
Let $V$ be an $n$-dimensional vector space and let $P \subset V$ be a compact convex polytope which is described as 
\begin{displaymath}
 P=\{x \in V: f_1(x) \geq 0,\ldots, f_m(x) \geq 0\},
\end{displaymath}
where $f_1,\ldots,f_m$ are affine functions. Let $d$ be the Finsler metric on $\inte P$ and endow the dual space $V^*$ with some norm. Then the map  
\begin{align*}
 \Phi: (\inte P,d) \to V^*\\
x \mapsto \sum_{i=1}^m \log f_i(x) df_i
\end{align*}
is a bilipschitz diffeomorphism. 
\end{Theorem}

Note that all norms on $V^*$ are equivalent, hence we may choose in the proof a Euclidean scalar product and identify $V^*$ with $V$. The map $\Phi$ can then be written as 
\begin{displaymath}
 \Phi(x)=\sum_{i=1}^m \log f_i(x) \grad f_i,
\end{displaymath}
where $\grad f_i \in V$ is the gradient of $f$. 

It is an interesting fact (which was communicated to me by J. Lagarias), that a similar map was studied in \cite{lagarias} in connection with trajectories for Karmarkar's linear programming algorithm.  

\subsection*{Acknowledgments}
We wish to thank Gautier Berck for pointing out several errors in an earlier version of this manuscript and Bruno Colbois and Patrick Verovic for useful discussions. 

\subsection*{Note} 
The main result of this manuscript, namely that the Hilbert geometry of a polytope is bilipschitz equivalent to a normed space, was also shown independently and with a different proof by Constantin Vernicos \cite{vernicos_lipschitz}.      
  

\section{Lipschitz continuity of $\Phi$}

Let us fix some notation. The Euclidean norm of a vector in $V$ will be denoted by $\|\cdot\|_2$, and $\|\cdot\|$ stands for the Finsler norm on $\inte P$. 

For $x \in \inte P$ and $0 \neq w \in T_x\inte P$, we have 
\begin{equation} \label{eq_norm_w}
 \|w\|=\frac12 \left(\frac{1}{t_1}+\frac{1}{t_2}\right),
\end{equation}
where $t_1,t_2 >0$ and $x+t_1 w, x-t_2 w \in \partial P$. 

We set $R:=\max_i \|\grad f_i\|_2$. By $D$ we will denote the (Euclidean) diameter of $P$. We set 
\begin{displaymath}
 P_i:=\{x \in V: f_i(x)=0\},
\end{displaymath}
which is an affine hyperplane.

\begin{Lemma} \label{lemma_equivalence_metrics}
\begin{enumerate}
\item For all $x,y \in \inte P$ we have 
\begin{displaymath}
 d(x,y) \geq \frac{\log 2}{D} \|x-y\|_2.
\end{displaymath}
In particular, for $x \in \inte P$ and $w \in T_x \inte P$ 
\begin{displaymath}
\|w\| \geq \frac{\log 2}{D} \|w\|_2.
\end{displaymath}
\item Let $K \subset \inte P$ be a compact set. There exists a constant $C_K>0$ such that for all $x \in K$ and all $w \in T_x\inte P$ 
\begin{displaymath}
	\|w\|_2 \geq C_K \|w\|.
\end{displaymath}
\end{enumerate}
\end{Lemma}

\proof
This follows easily from the definition of $d$ and from \eqref{eq_norm_w}. 
\endproof

\begin{Lemma} \label{lemma_lipschitz}
 The map $\Phi$ is Lipschitz continuous. 
\end{Lemma}

\proof
Let $x \neq y \in \inte P$ and $a_1,a_2 \in \partial P$ be the intersection points of the line through $x$ and $y$ with $\partial P$. Without loss of generality, let us assume that $a_1 \in P_1, a_2 \in P_2$ and that $x$ lies between $a_1$ and $y$. 

The Hilbert distance between $x$ and $y$ is given by 
\begin{align*}
 d(x,y) & =\frac12 \log \frac{\|y-a_1\|_2}{\|x-a_1\|_2} \frac{\|x-a_2\|_2}{\|y-a_2\|_2}\\
& =\frac12 \log \frac{f_1(y)}{f_1(x)}+\frac12 \log \frac{f_2(x)}{f_2(y)} \\
 & \geq \frac12 \max\left\{\left|\log \frac{f_1(x)}{f_1(y)}\right|,\left|\log \frac{f_2(x)}{f_2(y)}\right|\right\}.
\end{align*}

Since $a_1,a_2$ are the first intersection points with the boundary of $P$, for each $i \neq 1,2$ we have 
\begin{displaymath}
 \left|\log \frac{f_i(x)}{f_i(y)}\right| \leq \max\left\{\left|\log \frac{f_1(x)}{f_1(y)}\right|,\left|\log \frac{f_2(x)}{f_2(y)}\right|\right\}.
\end{displaymath}
It follows from these two inequalities that 
\begin{displaymath}
 d(x,y) \geq \frac12 \max\left\{\left|\log \frac{f_i(x)}{f_i(y)}\right|: i=1,\ldots,m\right\}.
\end{displaymath}

We now compute 
\begin{align*}
\| \Phi(x)-\Phi(y)\|_2 & = \left\|\sum_{i=1}^m \log f_i(x) \grad f_i- \sum_{i=1}^m \log f_i(y) \grad f_i\right\|_2\\
& \leq \sum_{i=1}^m \left|\log f_i(x)-\log f_i(y)\right| \cdot \|\grad f_i\|_2\\
& \leq 2mR d(x,y).
\end{align*}
\endproof

\begin{Lemma} $\Phi$ is injective.
\end{Lemma}

\proof
If $x \neq y$, there exists some $j$ with $f_j(x) \neq f_j(y)$. Noting that $\log$ is a strictly monotone function, we compute 
\begin{displaymath}
 \langle \Phi(x)-\Phi(y),x-y\rangle =\sum_{i=1}^m (\log f_i(x)-\log f_i(y)) (f_i(x)-f_i(y))> 0. 
\end{displaymath}

The injectivity of $\Phi$ follows. 
\endproof

\begin{Lemma} \label{lemma_norm_dphi}
 There exists a constant $C_1>0$ such that for all $x \in \inte P$ and all $w \in T_x\inte P$ 
\begin{displaymath}
 \|d\Phi(w)\|_2 \geq C_1 \|w\|_2.
\end{displaymath}
\end{Lemma}

\proof
The quadratic form $w \mapsto \sum_{i=1}^m \langle \grad f_i,w\rangle^2$ is positive definite, since the vectors $\grad f_i$ span $V$. Hence there is some constant $c_1$ with 
\begin{displaymath}
\sum_{i=1}^m \langle \grad f_i,w\rangle^2 \geq c_1 \|w\|_2^2. 
\end{displaymath}

Since $P$ is compact, there exists some real number $M>0$ with $f_i(x) \leq M$ for all $i$ and all $x \in \inte P$. Therefore,  
\begin{align*}
 \|d\Phi(w)\|_2 \cdot \|w\|_2 \geq \langle d\Phi(w),w\rangle & = \sum_{i=1}^m \frac{\langle \grad f_i,w\rangle^2}{f_i(x)} \geq \frac{c_1}{M} \|w\|_2^2.
\end{align*}
This proves the lemma, with $C_1=\frac{c_1}{M}$. 
\endproof

\section{Lipschitz continuity of $\Phi^{-1}$}

For a fixed polytope $P$, we consider two statements (A) and (B). 

\begin{itemize}
\item[(A)] There exists a constant $C>0$ such that for all $x \in \inte P$ and all $v \in T_x\inte P$ we have $\|d\Phi(v)\|_2 \geq C \|v\|$. 
\item[(B)] The map $\Phi:\inte P \to V$ is onto and bilipschitz. 
\end{itemize}

\begin{Proposition} \label{prop_a_implies_b}
If $P$ satisfies (A), then it also satisfies (B). 
\end{Proposition}

\proof
Set $W:=\Phi(\inte P) \neq \emptyset$. By Lemma \ref{lemma_norm_dphi} and the open mapping theorem, $W$ is open. Let $c:[0,1] \to W$ be a rectifiable curve and let $\tilde c:=\Phi^{-1} \circ c$ be its preimage under $\Phi$. Then 
\begin{equation} \label{eq_lipschitz_phi_inverse}
l(c) = \int_0^1 \|c'(s)\|_2 ds= \int_0^1 \|d\Phi(\tilde c'(s))\|_2 ds
\geq C \int_0^1 \|\tilde c'(s)\| ds = C l(\tilde c). 
\end{equation}  
If $W$ is not equal to $V$, there exists a curve $c:[0,1] \to W$ of finite length with $c(t) \in W$ for $0 \leq t < 1$ but $c(1) \notin W$. The preimage $\tilde c$ of $c|_{[0,1)}$ under $\Phi$ is a rectifiable curve of finite length in $(\inte P,d)$. Since this space is complete, $\tilde c$ can be extended to the whole interval $[0,1]$. By the Lipschitz property of $\Phi$, it follows that $c(1) = \Phi(\tilde c(1)) \in W$, a contradiction. Hence $\Phi$ is onto. 

Taking for $c$ the segment between two points, \eqref{eq_lipschitz_phi_inverse} implies that $\Phi^{-1}$ is Lipschitz with Lipschitz constant $\frac{1}{C}$.   
\endproof

\begin{Proposition} (A) is satisfied for {\it simple} polytopes. 
\end{Proposition}

\proof
We claim that there exists a number $\epsilon>0$ such that if we set, for $x \in \inte P$, 
\begin{displaymath}
I_\epsilon(x):=\{i \in \{1,\ldots,m\}: f_i(x)<\epsilon\}, 
\end{displaymath}
then 
\begin{displaymath}
 \bigcap_{i \in I_\epsilon(x)} P_i \neq \emptyset \quad \forall x \in \inte P. 
\end{displaymath}

If no such $\epsilon$ exists, we find a zero sequence $(\epsilon_l)$ and points $x_l \in \inte P$ with 
\begin{displaymath}
 \bigcap_{i \in I_{\epsilon_l}(x_l)} P_i = \emptyset.
\end{displaymath}
Passing to a subsequence, we may assume that the sets $I_{\epsilon_l}(x_l)$ are all equal to some $I$ and that the sequence $x_l$ converges to some point $x \in P$. Since $f_i(x_l) \to 0$, $x \in \cap_{i \in I} P_i$, which is a contradiction. 

For each $x \in \inte P$, the vectors $\{\grad f_i:i \in I_\epsilon(x)\}$ are linearily independent. Indeed, if $p$ is a vertex of the face $\cap_{i \in I_\epsilon(x)} P_i$, then the vectors $\{\grad f_j: f_j(p)=0\}$ span $V$ (see \cite{barv02}). On the other hand, since $P$ is simple, there are exactly $n$ such vectors and hence they are linearily independent.  

Let $C_P>0$ be such that whenever $\{\grad f_i,i \in I\}$ is some subset of \linebreak $\{\grad f_1,\ldots,\grad f_m\}$ of linearily independent vectors, then for all $\lambda_i \in \R$ we have 
\begin{displaymath}
 \left\|\sum_{i \in I} \lambda_i \grad f_i\right\|_2 \geq C_P \max_{i \in I} |\lambda_i|. 
\end{displaymath}
The existence of such a constant follows from the fact that any two norms on a finite-dimensional vector space are equivalent. 

Let $x \in \inte P$ and $w \in T_x \inte P$ with $\|w\|_2=1$. Let $t$ be the real number of minimal absolute value such that  
$x+tw \in \partial P$. Then $\|w\| \leq \frac{1}{|t|}$. 

We fix a number $\delta>0$ such that 
\begin{displaymath}
 \delta R \leq \epsilon, \frac{mR^2}{\epsilon}\leq \frac{C_P}{2\delta}. 
\end{displaymath}

Let us consider two cases. 
If $|t| \geq \delta$, then Lemma \eqref{lemma_norm_dphi} implies that 
\begin{displaymath}
 \|d\Phi(w)\|_2 \geq C_1 \|w\|_2=C_1 \geq C_1\delta \|w\|.
\end{displaymath}

If $|t|<\delta$, then $x+tw \in P_j$ for some $j$ and 
\begin{displaymath}
\frac{|\langle \grad f_j,w\rangle|}{f_j(x)}= \frac{1}{|t|}>\frac{1}{\delta}.
\end{displaymath}

It follows that $f_j(x)<\delta R \leq \epsilon$, i.e. $j \in I_\epsilon(x)$. 

We next compute that  
\begin{align*}
 \|d\Phi(w)\|_2 & =\left\|\sum_{i=1}^m \frac{\langle \grad f_i,w\rangle}{f_i(x)} \grad f_i\right\|_2\\
& \geq \left\|\sum_{i \in I_\epsilon(x)} \frac{\langle \grad f_i,w\rangle}{f_i(x)} \grad f_i\right\|_2-\left\|\sum_{i \not\in I_\epsilon(x)} \frac{\langle \grad f_i,w\rangle}{f_i(x)}\grad f_i\right\|_2\\
& \geq C_P \frac{|\langle \grad f_j,w\rangle|}{f_j(x)}-\frac{m R^2}{\epsilon}\\
& \geq \frac{C_P}{2} \frac{|\langle \grad f_j,w\rangle|}{f_j(x)}\\
& \geq  \frac{C_P}{2} \|w\|.
\end{align*}

Property (A) now follows with $C:=\min\left\{C_1\delta,\frac{C_P}{2}\right\}$. 
\endproof

Noting that every plane convex polygon is simple, we obtain that Theorem \ref{thm_main_thm_exact} holds true in dimension $2$. We now proceed by induction on the dimension of $P$. We first prove property (A) for polyhedral cones. 

In the definition of the Hilbert metric, we supposed $K$ to be compact and convex. However, even if $K$ is only closed and convex and does not contain any straight line, then the definition makes sense (in this case, one of $a_1$ or $a_2$ may be at infinity). 

\begin{Proposition} \label{prop_cones}
Let $P \subset V$ be a polyhedral cone of the form 
\begin{displaymath}
P=\{x \in V:f_1(x) \geq 0,\ldots,f_m(x) \geq 0\}
\end{displaymath}
where $f_1,\ldots,f_m$ are linear functions on $V$. Suppose that $P$ does not contain any line, but has non-empty interior. Define a map 
\begin{align*}
 \Phi: \inte P & \to V\\
 x & \mapsto \sum_{i=1}^m \log f_i(x) \grad f_i.
\end{align*}
Then there exists a constant $C>0$ such that for each $x \in \inte P$ and each $w \in T_x \inte P$, we have 
\begin{equation} \label{eq_lipschitz_cone}
\|d\Phi(w)\|_2 \geq C \|w\|.
\end{equation} 
\end{Proposition}

\proof
Let $u:=\sum_{i=1}^m \grad f_i$ and $E_0:=u^\perp$. Let $x_0 \in \inte P$ and set $E:=x_0 + E_0$. Then $P^E:=P \cap E$ is a compact polytope of dimension $n-1$.  

By an easy homogeneity argument, it suffices to prove \eqref{eq_lipschitz_cone} for $x \in \inte P^E$. We let $f_i^E$ denote the restriction of $f_i$ to $E$. With $\pi:V \to E_0$ being the orthogonal projection, we have $\grad f_i^E=\pi(\grad f_i)$. 

Since $P^E:=P \cap E$ is of dimension $n-1$, there is a constant $C^E$ such that for all $x \in \inte P^E$ and all $w \in T_x \inte P^E$ we have 
\begin{equation} \label{eq_bilipschitz_restriction}
 \left\|\sum_{i=1}^m \frac{\langle \grad f_i,w\rangle}{f_i(x)} \grad f_i^E\right\|_2 \geq C^E \|w\|.
\end{equation}

Fix a positive constant $c_1$ with 
\begin{displaymath}
(1-c_1) \|u\|_2 -mR c_1 >0.
\end{displaymath} 

Let $x \in \inte P^E$ and let $w \in T_x\inte P$ be of Euclidean norm $1$. We write $w=w_1+w_2$ with $w_1$ parallel to $E$ and $w_2$ parallel to the line $\R \cdot x$, say $w_2=\rho x$. Note that $\|w_2\| =\frac12 |\rho|$. 

By \eqref{eq_bilipschitz_restriction},  
\begin{displaymath}
\|\pi(d\Phi(w_1))\| = \left\|\sum_{i=1}^m \frac{\langle \grad f_i,w_1\rangle}{f_i(x)} \grad f_i^E\right\|_2 \geq C^E\|w_1\|. 
\end{displaymath}

Next we compute that 
\begin{displaymath}
d\Phi(w_2) = \sum_{i=1}^m \frac{\langle \grad f_i,w_2\rangle}{f_i(x)}\grad f_i= \rho u.
\end{displaymath}

In particular
\begin{displaymath} 
 \pi(d\Phi(w_2))=0; \quad  \|d\Phi(w_2)\|_2 =2\|w_2\| \cdot \|u\|_2.
\end{displaymath}

If $\|w_1\| \geq c_1 \|w\|$, then we obtain 
\begin{displaymath}
\|d\Phi(w)\|_2 \geq \|\pi(d\Phi(w))\|_2=\|\pi(d\Phi(w_1))\|_2 \geq C^E \|w_1\| \geq c_1 C^E \|w\|. 
\end{displaymath}

If $\|w_1\| < c_1 \|w\|$, then by triangle inequality $\|w_2\| \geq (1-c_1) \|w\|$ and, using that $\Phi$ is Lipschitz with Lipschitz constant $2mR$, where $R=\max_i \|\grad f_i\|_2$ (see Lemma \ref{lemma_lipschitz}), we get   
\begin{align*}
\|d\Phi(w)\|_2 & \geq \|d\Phi(w_2)\|_2-\|d\Phi(w_1)\|_2\\
& \geq 2\|w_2\| \cdot \|u\|_2-2mR \|w_1\|\\
& \geq 2(1-c_1) \|u\|_2 \cdot \|w\|-2mR c_1 \|w\|.
\end{align*} 
Thus \eqref{eq_lipschitz_cone} is satisfied with $C:=\min\{c_1C^E,2(1-c_1) \|u\|_2 -2mR c_1\}$. 
\endproof

\begin{Proposition} \label{lemma_property_a_polytopes}
Property (A) is satisfied for each polytope $P$ of dimension $n$.
\end{Proposition}

\proof
Let us prove by induction on $k=0,\ldots,n-1$ the following statement:
\begin{itemize}
\item[(A$_k$)] For every $k$-dimensional face $F$ of $P$, there exists an open neighborhood $U$ in $V$ and a constant $C_F$ such that for all $x \in \inte P \cap U$ and all $v \in T_x\inte P$ we have $\|d\Phi(v)\|_2 \geq C_F \|v\|$. 
\end{itemize}

It is clear that (A$_{n-1}$), Lemma \ref{lemma_equivalence_metrics} and Lemma \ref{lemma_norm_dphi} imply (A). The induction start will be the empty case $k=-1$. 

Suppose now $k\geq 0$. Let $F=\cap_{i \in I} P_i$ be a $k$-face of $P$. By induction hypothesis, we may assume that there is an open neighborhood $U'$ of the $k-1$-skeleton of $P$ and a constant $C_2$ such that for all $x \in \inte P \cap U'$ and all $v \in T_x\inte P$ we have $\|d\Phi(v)\|_2 \geq C_2 \|v\|$.

On the compact set $F \setminus U'$, the continuous functions $f_j, j \notin I$ are strictly positive and hence strictly larger than some constant $\tau>0$. Set  
\begin{displaymath}
	U'':=\{x \in V: f_j(x)>
	\tau \quad \forall j \notin I\}.
\end{displaymath}
Then $U''$ is open and $U' \cup U''$ is an open neighborhood of $F$.

Let $F_0$ be the $k$-dimensional linear space parallel to $F$  and set $\bar V:=V/F_0$. The affine functions $f_i:V \to \R$ induce linear functions $\bar f_i:\bar V \to \R$. Define a polyhedral cone 
\begin{displaymath}
 \bar P:=\left\{\bar x \in \bar V: \bar f_i(\bar x) \geq 0, i \in I\right\}. 
\end{displaymath}
 
By Proposition \ref{prop_cones}, applied to $\bar P$, there exists a constant $C_3>0$ such that for all $x \in \inte P \cap U''$ and all $u \in T_x\inte P \cap F_0^\perp$ we have  
\begin{equation} \label{eq_bound_conical_part}
\left\|\sum_{i \in I} \frac{\langle u,\grad f_i\rangle}{f_i(x)}\grad f_i\right\|_2 \geq C_3 \overline{\|u\|}.
\end{equation}
Here $\overline{\|u\|}$ denotes the Finsler norm with respect to the polyhedral cone $P'$. 

Let $x \in \inte P \cap U''$ and $w \in T_x \inte P$ with $\|w\|_2=1$. We write $w=w_1+w_2$ with $w_1 \in F_0^\perp$ and $w_2 \in F_0$. Then 
\begin{displaymath}
 \|d\Phi(w_2)\|_2=\left\|\sum_{j \notin I} \frac{\langle w_2,\grad f_j\rangle}{f_j(x)}\grad f_j\right\|_2 \leq \frac{m R^2}{\tau}.
\end{displaymath}

Fix a sufficiently large positive constant $c_1$ with $\frac{R}{c_1\tau-R} < 1$ and
\begin{displaymath}
 C_4:=\frac{C_3}{2} \left(1-\frac{R}{\tau c_1}\right)-\frac{2mR^2}{\tau c_1}>0. 
\end{displaymath}

We consider two cases.

Case 1: $\|w\| \leq c_1$. Then $\|d\Phi(w)\|_2 \geq C_1\|w\|_2 \geq \frac{C_1}{c_1} \|w\|$.
 
Case 2: $\|w\| \geq c_1$. Let $t$ be of minimal absolute value such that $x+tw_2 \in \partial P$, say $x+tw_2 \in P_j$. Since $w_2$ is parallel to $F$, we have $j \notin I$. From
\begin{displaymath}
 \tau < f_j(x) =|f_j(x+tw_2)-f_j(x)|=|t| \cdot |\langle \grad f_j,w_2\rangle| \leq R|t|
\end{displaymath}
we deduce that $\|w_2\| \leq \frac{R}{\tau}$. The triangle inequality now yields 
\begin{displaymath}
 \|w_1\| \geq \|w\|-\|w_2\| \geq \|w\|-\frac{R}{\tau c_1}\|w\| \geq \frac{c_1 \tau -R}{\tau}. 
\end{displaymath}

Let $s$ be of minimal absolute value with $x+sw_1 \in \partial P$, say $x+sw \in P_i$. Then $|s| \leq \frac{1}{\|w_1\|} \leq \frac{\tau}{c_1 \tau -R}$ and therefore $f_i(x) \leq  R|s| \leq \frac{R\tau}{c_1\tau-R}<\tau$, hence $i \in I$. 

By \eqref{eq_norm_w}, we have $\|w_1\| \leq \frac{1}{|s|}$. The Finsler norm with respect to the cone $\bar P$ satisfies $\overline{\|w_1\|} \geq \frac{1}{2|s|}$. By \eqref{eq_bound_conical_part} we get 
\begin{displaymath}
 \left\|\sum_{i \in I} \frac{\langle w_1,\grad f_i\rangle}{f_i(x)}\grad f_i\right\|_2 \geq \frac{C_3}{2} \|w_1\|. 
\end{displaymath}

Finally, we obtain 
\begin{align*}
 \|d\Phi(w)\|_2 & \geq \|d\Phi(w_1)\|_2-\|d\Phi(w_2)\|_2\\
& \geq \left\|\sum_{i \in I} \frac{\langle w_1,\grad f_i\rangle}{f_i(x)}\grad f_i\right\|_2-\left\|\sum_{j \notin I} \frac{\langle w_1,\grad f_j\rangle}{f_j(x)}\grad f_j\right\|_2\\
& \quad -\frac{mR^2}{\tau}\\
& \geq \frac{C_3}{2} \|w_1\|-\frac{2mR^2}{\tau}\\
& \geq \frac{C_3}{2} \left(1-\frac{R}{\tau c_1}\right) \|w\|-\frac{2mR^2}{\tau}\\
& = C_4 \|w\|+\frac{2mR^2}{\tau c_1} \|w\|-\frac{2mR^2}{\tau}\\
& \geq C_4 \|w\|.
\end{align*}

We infer that A$_k$ holds true with 
\begin{displaymath}
 U:=U' \cup U''; \quad C_F:=\min\left\{C_2,\frac{C_1}{c_1},C_4\right\}.
\end{displaymath}
\endproof

Theorem \ref{thm_main_thm_exact} clearly follows from Propositions \ref{lemma_property_a_polytopes} and \ref{prop_a_implies_b}.

\end{document}